\documentclass[leqno, draft]{amsart}
\usepackage[c5paper, text={126truemm,178truemm},centering]{geometry}
\usepackage{xcolor}
\usepackage{tikz}
\usepackage{comment}
\usepackage{pict2e}
\usepackage[utf8]{inputenc}
\usetikzlibrary{calc,arrows,decorations.markings,intersections,matrix}
\newcommand\includetikz[2][x=1mm,y=1mm]{%
 \IfFileExists{#2.tikz}{%
	\begin{tikzpicture}[#1]\node at (0,0) {\input{#2.tikz}};\end{tikzpicture}%
 }{%
	\@latex@error{No usable file `#2.tikz' can be found}%
	 {I could not locate the file...}%
 }
}
\def\Bf#1{\ifmmode\boldsymbol{#1}\else{\rmfamily\bfseries#1}\fi}
\def\<{\langle}
\def\>{\rangle}

\usepackage[backref=page]{hyperref}
\hypersetup{colorlinks, citecolor=[rgb]{0.0,0.3,0.1}, linkcolor=[rgb]{0.3,0.0,0.1},
						urlcolor=[rgb]{0.0,0.0,0.4}}
\makeatletter
\def\BR@backref#1{{\upshape
 \begingroup
  \csname @safe@activestrue\endcsname
  \expandafter\providecommand\csname brc@#1\endcsname{0}%
  \expandafter\providecommand\csname brcd@#1\endcsname{0}%
  \csname @safe@activesfalse\expandafter\endcsname
  \ifBR@BackrefAlt
   \ifx\backrefentrycount\BR@BackrefEntryCountUnused
   \else\BR@PopulateEntryCount{#1}\fi
   \expandafter\backrefalt\csname brc@#1\expandafter\endcsname
   \csname brl@#1\expandafter\endcsname
   \csname brcd@#1\expandafter\endcsname
   \csname brld@#1\endcsname
  \else
   \expandafter\backref\csname br@#1\expandafter\endcsname
  \fi
 \endgroup
 \par
}}
\makeatother
\renewcommand*{\backref}[1]{\ifx#1\relax\relax\else{\rm~$\langle$#1$\rangle$}\fi}

\theoremstyle{plain}
\newtheorem{theorem}{Theorem}[section]

\theoremstyle{definition}

\newtheorem{remark}[theorem]{Remark}

\numberwithin{figure}{section}

\title[a Plancherel-type formula]{                   
					Plancherel formula for the attenuated Radon transform
}

\author[E. Dinnyés, T. Ódor]{Enikő Dinnyés${}^{*,\#}$ \and Tibor Ódor${}^{**,\#}$}
\thanks{${}^*$ Department of Probability and Statistics, Eötvös Loránd University, Budapest}
\thanks{${}^{**}$ Department of Geometry, Bolyai Institute, University of Szeged}
\thanks{${}^{\#}$ Rényi Institute of Mathematics, Budapest}
\thanks{This work on the project leading to this application has received
funding from the European Research Council (ERC) under the European Union’s Horizon 2020 research and innovation programme (grant agreement No. 741420).}

\begin{document}


\begin{abstract}
Motivated by stereology, based on Novikov's inversion formula, we prove a Plancherel-type formula for the attenuated 
Radon transform.
\end{abstract}

\maketitle

\section{\label{sec:intro}
                              Introduction
}

In this article we construct a Plancherel-type formula for the attenuated Radon transform. We shall compute $\int_{\mathbb R^2} f(x)\, g(x)\, dx$ from the attenuated Radon transform $R_af$, for any $g$ of our choice (known exactly or with arbitrary precision).

The problem is motivated by methods in stereology~\cite{Baddeley}.
Using classical Radon transform~\cite{Helgason}, it is easy to estimate the area or volume of a body, or an integral of a function with compact support, as the expected value of sections, or integrals over lines or planes, just by the Fubini theorem. But this does not work for the attenuated Radon transform due to its complicated weights over the line we integrate on.

For deriving a Plancherel-type formula for the attenuated Radon transform, we use Natterer's version~\cite{Natterer} of Novikov's inversion formula~\cite{Novikov}.

\bigskip
Let $a\colon \mathbb R^2\to \mathbb R$ be a sufficiently smooth and sufficiently decaying function (as in \cite{Novikov}), e.g. an integrable $C^1$ function whose partial derivatives are also integrable. For $\omega \in \mathbb S^1$ and $x\in \mathbb R^2$ we define
\begin{equation}
(Da)(x, \omega) = \int^\infty_0 a(x + t\,\omega) \, dt.
\end{equation}
We will use $(Da)(x, \omega^\perp)$, with $\omega^\perp$ defined below. 

\begin{remark} Depending on the context, the symbol $\omega$ can be understood as the angle of the vector $\omega$ with the first coordinate axis, or the vector itself, that is, $\omega=(\cos\omega, \sin\omega)$. There is a one-to-one correspondance between the vector $\omega\in\mathbb S^1$ and $\omega\in[0,2\pi)$ (geometrically equivalent to an angle). When the argument of the function to be differentiated is $x\in\mathbb R^2$, then $\partial_\omega$ (or $\partial_{\underline\omega}$) means directional derivative in the direction $\omega$. When one of the arguments of the function to be differentiated is $\omega$, then $\partial_\omega$ means partial derivative w.r.t. the angle $\omega$. Occasionally we underline the name of vectors for easier understanding.
\end{remark}

The \emph{attenuated Radon transform} $R_a: [0,2\pi) \times \mathbb{R} \rightarrow \mathbb{R}$ (or equivalently, $R_a: (-\pi,\pi] \times \mathbb{R} \rightarrow \mathbb{R}$) is defined as 
\begin{equation}
(R_a f)(\omega, p) = \int_{\< x, \omega\> = p} e^{-Da( x,\, \omega^\perp)}\, f( x)\, dx,
\end{equation}
where $dx$ is the Lebesgue measure on the line $l(\omega, p) := \left\{ x \in\mathbb{R}^2 : \< x,\omega\> = p \right\}$, and 
$\omega^\perp$ is 
$(\cos \omega, \sin\omega)^\perp = (\cos (\omega + \frac{\pi}{2}), \sin(\omega + \frac{\pi}{2})) = (-\sin\omega, \cos\omega)$. 
With $a=0$ we get the classical Radon transform of $f$.

We can rewrite this definition as 
\begin{equation}
(R_a f)(\omega, p) = \int^\infty_{-\infty} e^{-Da(p\, \omega + u\, \omega^\perp,\, \omega^\perp)}\, f(p\, \omega + u\, \omega^\perp)\, du.
\end{equation}

Note that 
the line $l(\omega, p)$ is the same set as $l(-\omega, -p)$, but admits different orientation. Also note that although defined on the same set, the integral $(R_a f)(\omega, p)$ usually differs from $(R_a f)(-\omega, -p)$, unlike the special case of the classical Radon transform where $Rf(\omega,p)=Rf(-\omega,-p)$.

\bigskip
\section{\label{sec:plancherel}
                           								The Plancherel formula 
}

Let us define the function $h(\omega,p) = \frac{1}{2}(I + i \mathcal{H})Ra(\omega,p)$, where $R a(\omega,p)$ is the classical Radon transform of the function $a$ along the line $l(\omega,p)$, and $\mathcal{H}$ is the Hilbert transform, defined by $\mathcal{H}g(p) = \frac{1}{\pi} \int^\infty_{-\infty} \frac{g(t)}{p-t}\, dt$ where the integral is understood as a Cauchy principal value. 
We have the following reconstruction formula due to Natterer \cite{Natterer} based on Novikov's \cite{Novikov}, if $a$ and $f$ are sufficiently smooth functions decaying sufficiently fast at infinity:

\begin{equation}
\begin{split}
f(x) = \frac{1}{4\pi} {\rm div}\, \int_{\mathbb S^1} \underline\omega \,\left\{ e^{Da(x, \omega^\perp)} ({\rm Re}\; e^{-h} \mathcal{H} e^h R_af)(\omega, \< x, \omega\>)\right\} \, d\omega = \\
= \frac{1}{4\pi} \frac{\partial}{\partial x_1} \int_{\mathbb S^1} \cos\omega \,\left\{ e^{Da(x, \omega^\perp)} ({\rm Re}\; e^{-h} \mathcal{H} e^h R_af)(\omega, \< x, \omega\>)\right\} \, d\omega \,+ \\
+ \frac{1}{4\pi} \frac{\partial}{\partial x_2}\int_{\mathbb S^1} \sin\omega \,\left\{ e^{Da(x, \omega^\perp)} ({\rm Re}\; e^{-h} \mathcal{H} e^h R_af)(\omega, \< x, \omega\>)\right\} \, d\omega = \\
= \frac{1}{4\pi} \int_{\mathbb S^1} \partial_{\underline\omega}(e^{Da(x,\omega^\perp)})({\rm Re}\; e^{-h} \mathcal{H} e^h R_af)(\omega, \< x, \omega\>) \,d\omega\,+\\
+\frac{1}{4\pi} \int_{\mathbb S^1} e^{Da(x,\omega^\perp)}\; \partial_p ({\rm Re}\; e^{-h} \mathcal{H} e^h R_af)(\omega, \< x, \omega\>) \, d\omega,
\end{split}
\end{equation}
as ${\rm div} (\underline\omega\, c(\underline x)) = \cos\omega \, \frac{\partial}{\partial x_1}c(\underline x) + \sin\omega\, \frac{\partial}{\partial x_2}c(\underline x)$;
and $\partial_{\underline\omega} = \cos\omega\;\frac{\partial}{\partial x_1} + \sin\omega\;\frac{\partial}{\partial x_2}$ is the directional derivative of a function (defined in $\mathbb{R}^2$)  in the direction $\omega$, which is now acting on the first variable of $e^{Da(x,\omega^\perp)}$; 
and $\partial_p$ is the partial derivative of a function defined on $\mathbb{P}^2=\mathbb{S}^1\times\mathbb{R}_+$, the space of all straight lines, with respect to its second variable (the distance from the origin). The latter comes from $\frac{\partial}{\partial x_1}$ and $\frac{\partial}{\partial x_2}$ acting on the scalar product 
$\< x, \omega\>$ as an inside function, with derivatives $\cos\omega$ and $\sin\omega$, respectively, multiplying $\cos\omega$ and $\sin\omega$ that were already there, so adding up to 1.

Note  that if $p$ is changing with the direction $\omega$ fixed, that has the same geometric meaning as the directional derivative in the direction $\omega$.

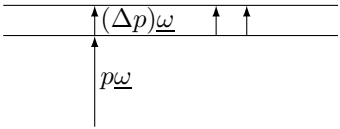
\setlength{\unitlength}{0.4 cm}
\begin{picture}(12,6)(-3.5,-3.5)
\put(-3,0){\line(1,0){11}}
\put(-3,1){\line(1,0){11}}
\put(0,0){\vector(0,1){1}}
\put(0.2, 0.5){\makebox(0,0)[l]{$(\Delta p)\underline{\omega}$}}
\put(4,0){\vector(0,1){1}}
\put(5,0){\vector(0,1){1}}
\put(0,-3){\vector(0,1){3}}
\put(0.2,-1.5){\makebox(0,0)[l]{$p\underline{\omega}$}}
\end{picture}

Our goal is to compute $\int_{\mathbb R^2} f(x)\, g(x)\, dx$ from the attenuated Radon transform $R_af$, knowing $g$ exactly (with arbitrary precision). An application of this approach is when $g$ is the indicator function of a small area. Then knowing the integral of $f\cdot g$ is close to reconstruction of $f$ locally.

Now we multiply the above expression of $f(x)$ by $g(x)$ and integrate it on $\mathbb{R}^2$, then change the order of integration. The differentiation with respect to $p$ when $u$ and $\omega$ are fixed, and the directional derivative in the direction $\omega$, are exactly the same, so the two terms of our previous expression for $f(x)$ can be contracted to a total derivative (using the substitutions 
$\< x,\omega\> =p$ and $\< x,\omega^\perp\> =u$):
\begin{equation}
\begin{split}
\int_{\mathbb{R}^2} f(x)\, g(x)\, dx =\\
\int_{\mathbb{R}^2} \bigg\{ \frac{1}{4\pi} \int_{\mathbb S^1}\partial_{\underline\omega} e^{Da(x,\omega^\perp)}\,({\rm Re}\; e^{-h} \mathcal{H} e^h R_af)(\omega, \< x, \omega\>)+ \\
+ e^{Da(x,\omega^\perp)}\; \partial_p ({\rm Re}\; e^{-h} \mathcal{H} e^h R_af)(\omega, \< x, \omega\>) \, d\omega \bigg\}\, g(x)\,dx =\\
= \frac{1}{4\pi} \int_{\mathbb S^1} \int_{-\infty}^\infty \int_{-\infty}^\infty \bigg\{ \partial_p\{ e^{\int_0^\infty a(p\omega+u\omega^\perp + t\omega^\perp)dt}\}\,
({\rm Re}\; e^{-h} \mathcal{H} e^h R_af)(\omega, p)\,
+\\
+
e^{\int_0^\infty a(p\omega+u\omega^\perp + t\omega^\perp)dt}\, \partial_p ({\rm Re}\; e^{-h} \mathcal{H} e^h R_af)(\omega, p)\bigg\} \, g(p\omega+u\omega^\perp) \,dp\, du\, d\omega =\\
\frac{1}{4\pi} \int_{\mathbb S^1} \int_{-\infty}^\infty \int_{-\infty}^\infty \partial_p\left\{e^{\int_0^\infty a(p\omega+u\omega^\perp + t\omega^\perp)dt}\;{\rm Re}\; e^{-h} \mathcal{H} e^h R_af(\omega, p)\right\}\cdot \\
\cdot \; g(p\omega+u\omega^\perp)\,dp\, du\,d\omega\,=\\
 -\frac{1}{4\pi} \int_{\mathbb S^1} \int_{-\infty}^\infty \int_{-\infty}^\infty \left\{e^{\int_0^\infty a(p\omega+u\omega^\perp + t\omega^\perp)dt}\;{\rm Re}\; e^{-h} \mathcal{H} e^h R_af(\omega, p)\right\}\cdot \\
 \cdot\; \partial_p\left\{g(p\omega+u\omega^\perp)\right\}\,dp\, du\,d\omega\,=\\
\end{split}
\end{equation}
\begin{equation}
\begin{split}
= -\frac{1}{4\pi} \int_{\mathbb S^1} \int_{-\infty}^\infty {\rm Re}\; e^{-h} \mathcal{H} e^h R_a f(\omega, p)\;\cdot\\
 \cdot \left\{ \int_{-\infty}^\infty e^{\int_0^\infty a(p\omega+u\omega^\perp + t\omega^\perp)dt}\,(\partial_{\underline\omega}g) (p\omega+u\omega^\perp)\,du\right\} dp\,d\omega\,=\\
 -\frac{1}{4\pi} \int_{\mathbb S^1} \int_{-\infty}^\infty {\rm Re}\; e^{-h} \mathcal{H}\left\{ e^h R_af(\omega, p)\right\}\cdot R_a(\partial_{\underline\omega}g)(\omega, p)\,dp\,d\omega\,=\\
\frac{1}{4\pi} \int_{\mathbb S^1} \int_{-\infty}^\infty {\rm Re}\; e^{h} R_a f(\omega, p) \; \mathcal{H} \left\{e^{-h} R_a(\partial_{\underline\omega}g)(\omega, p)\right\} dp\,d\omega.\\
\end{split}
\end{equation}
Here we used integration by parts w.r.t. $p$, and the following observation about the Hilbert transform:
\begin{equation}
\begin{split}
\int f(x)\mathcal{H}g(x)dx = \int f(x) \left\{\int \frac{g(r)}{x-r}dr\right\} dx = \int \int \frac{f(x)g(r)}{x-r}dr\,dx =\\
-\int \int \frac{g(r)f(x)}{r-x}dx\,dr = -\int g(r)\left\{\int \frac{f(x)}{r-x}dx\right\} dr = -\int g(r)\mathcal{H}f(r)dr.\\
\end{split}
\end{equation}
This way both the differentiation and the Hilbert transform was "transferred" from the unknown $f$ to the known $g$.
We have proven the following theorem, with the conditions of Novikov's inversion formula.

\begin{theorem}
The Plancherel formula for the complete attenuated Radon transform $R_a$ is as follows:
\begin{equation}
\begin{split}
\int_{\mathbb R^2} f(x)\, g(x)\, dx = \\ \frac{1}{4\pi}\, 
\int_{\mathbb S^1}\int^\infty_{-\infty} (R_a f)(\omega, p) \, {\rm Re}\left[ e^{h(\omega, p)} 
\mathcal{H}_p \left\{ e^{-h(\omega, p)} \, (R_a \partial_\omega g)(\omega, p)\right\} \right] dp\, d\omega. \;\;\;\; \square
\end{split}
\end{equation}
\end{theorem}  

If $g$ is known to arbitrary precision, because for example it is of our choice, e.g. it can be constant on a big circular disc containing the support of $f$, then this computes the integral of $f$, which is of major concern in stereology. 
If $a = 0$, we deal with the classical Radon transform, and through the Fubini theorem, having the integral of $f$ is trivial. But not so, if $a$ is non-zero. 
To our knowledge, this is the first formula for that problem. Observe that we need no derivatives or Hilbert transform of $R_a f$. 
If $g$ is known, then $ R_a (\partial_\omega g)$ can be computed. So our formula is numerically stable and e.g. Monte Carlo and related methods can be applied.

Although our formula is "ideal" for the needs of stereology, for the same reason, due to its asymmetry in $f$ and $g$ it is not so for theoretical reasons, like extension of the attenuated Radon transform to $L^2$ spaces. Finding a useful symmetric version (except the obvious one, when we express $g$ from $R_a g$ using Novikov's formula) is still open.

\end{document}